\documentclass[11pt]{article}
\usepackage{amsthm}
\usepackage{latexsym,amssymb,amsmath,enumerate}
\usepackage{graphicx,float,color,fancybox,shapepar,setspace,hyperref}
\usepackage{subfigure}
\usepackage{pgf,tikz}
\usetikzlibrary{arrows}
\newtheorem{thm}{Theorem}[section]
\newtheorem{conj}[thm]{Conjecture}

\newtheorem{claim}[thm]{Claim}
\newtheorem{lem}[thm]{Lemma}
 
 \def\dfn#1{{\sl #1}}
 \def\qed{\hfill\square}

\def\qed{ \hfill $\blacksquare$}
\def\pf{\noindent{\emph{Proof}.}~~}
\def\less{\setminus}
\def\es{\emptyset}
\usepackage[affil-it]{authblk}

\begin{document}

\title{Gallai-Ramsey number of odd cycles with chords}

\author{Fangfang  Zhang$^{1,}$\thanks{This work was done while the first author visited the University of Central Florida as a visiting student. The hospitality of the
hosting institution is greatly acknowledged.    The visit was funded by the Chinese Scholarship Council.  E-mail address: fangfangzh@smail.nju.edu.cn.}  \,\,,  Zi-Xia Song$^{2,}$\thanks{Partially supported by the National   Science  Foundation of China under Grant No. DMS-1854903.   E-mail address:  Zixia.Song@ucf.edu. }\,\,,    Yaojun Chen$^{1,}$\thanks{Yaojun Chen and Fangfang Zhang are partially supported by the National Natural Science Foundation
of China  under grant numbers 11671198 and 11871270. E-mail address:  yaojunc@nju.edu.cn.}}
 \affil{ { \small {$^1$Department of Mathematics, Nanjing University, Nanjing 210093, China}}\\
  { \small {$^2$Department  of Mathematics, University of Central Florida, Orlando, FL 32816, USA}}
 }

\date{}
\maketitle

\begin{abstract}
   A Gallai coloring is a coloring of the edges of a complete graph without
rainbow triangles, and a    Gallai $k$-coloring is a Gallai coloring that uses at most $k$ colors.   For an integer $k\geq 1$,   the  Gallai-Ramsey number $GR_k(H)$ of a given graph $H$ is the least positive  integer $N$ such that every Gallai $k$-coloring of the complete graph $K_N$   contains a monochromatic copy of $H$.   Let $C_m$ denote the cycle  on  $m\ge4$ vertices        and    let $\Theta_m$ denote  the family of   graphs obtained from   $C_m$ by adding  an additional edge joining two non-consecutive vertices.  We prove that  $GR_k(\Theta_{2n+1})=n\cdot 2^k+1$ for all $k\geq 1$ and $n\geq 3$. This implies that 
 $GR_k(C_{2n+1})=n\cdot 2^k+1$  all $k\geq 1$ and $n\geq 3$. Our result  yields   a unified proof for the Gallai-Ramsey number of  all odd cycles on at least five vertices.   \medskip
 
\noindent{\bf Keywords}: Gallai coloring,   Ramsey theory,  cycles, rainbow triangle   

\noindent{\bf 2020 Mathematics Subject Classification}: 05C55;  05D10; 05C15
\end{abstract}
 
\section{Introduction}
 In this paper we consider graphs that are finite, simple and undirected.      We use $P_m$,  $C_m$ and $K_m$ to denote the path,    cycle and  complete graph  on $m$ vertices, respectively. For $m\ge4$, let $\Theta_m$ denote  the family of   graphs obtained from   $C_m$ by adding  an additional edge joining two non-consecutive vertices.  
For any positive integer $k$, we write  $[k]$ for the set $\{1,2, \ldots, k\}$. 
 Given an integer $k \ge 1$ and graphs $H_1,   \ldots, H_k$, the classical \dfn{Ramsey number} $R(H_1,    \ldots, H_k)$   is  the least    integer $N$ such that every $k$-coloring of  the edges of      $K_N$  contains  a monochromatic copy of  $H_i$ in color $i$ for some $i \in [k]$.  When $H = H_1 = \dots = H_k$, we simply write   $R_k(H)$ to denote the \dfn{$k$-color Ramsey number} of $H$.    In the seminal paper of Ramsey \cite{Ramsey}, it is shown that Ramsey numbers are finite. This was   rediscovered independently by Erd\H{o}s and Szekeres  \cite{ES}.  Since the 1970s, Ramsey theory has grown into one of the most  active areas of research   in combinatorics, overlapping variously with graph theory, number theory, geometry and logic. However, determining Ramsey numbers is  notoriously difficult   in general. \medskip
 
  Let $(G,\tau)$ denote a $k$-colored  complete graph, where $G$ is a complete graph and $\tau: E(G)\rightarrow [k]$. We say $(G,\tau)$ is \dfn{$\mathcal{F}$-free} if $G$ does not contain a monochromatic copy of a   graph in a given family $\mathcal{F}$ under  the $k$-coloring $\tau$; when $\mathcal{F}=\{F\}$, we simply say $(G,\tau)$ is \dfn{$F$-free}. By abusing notation, we say  $(G,\tau)$ contains a \dfn{monochromatic copy of   $\mathcal{F}$} if $G$  contains a monochromatic copy of a   graph in  $\mathcal{F}$ under $\tau$.  One of the earliest and well-known problems is that of  determining the Ramsey number    $R_k(C_m)$.  For odd cycles $C_{2n+1}$,  Erd\H{o}s and Graham~\cite{EG} observed that one can naturally construct a $C_{2n+1}$-free $(G, \tau)$  by induction. Indeed  when $k = 1$,  simply take $(G, \tau)$ to be   a $1$-colored  $K_{2n}$; for $k > 1$,  let $(G, \tau)$ be obtained by taking two disjoint copies of the construction for $k -1$ and color every edge between the two copies with a new color. This construction showed that 
 $R_k(C_{2n+1})\ge n\cdot 2^k+1 $ for all $k\ge1$ and $n\ge2$, which 
  led Bondy and Erd\H{o}s \cite{BE} to make the following conjecture.

\begin{conj}[\cite{BE}]\label{OddCycle} For all  $k\ge1$ and $n\ge2$,
$R_k(C_{2n+1})=n\cdot 2^k+1. $
 
 \end{conj}

 Recently, Jenssen and Skokan~\cite{JenssenSkokan}  showed that Conjecture~\ref{OddCycle} is true for all fixed $k$ and all $n$ sufficiently large. 
 However,   for all fixed $n$ and all $k$ sufficiently large,  Day and Johnson~\cite{DayJohnson} recently  showed that  $R_k(C_{2n+1}) > 2n\cdot (2+\epsilon)^{k-1}$ for some constant $\epsilon= \epsilon(n)> 0$, and so Conjecture~\ref{OddCycle} is false when $n$ is small with respect to $k$. 
For further results in this direction, we refer the reader to Graham, Rothchild and Spencer~\cite{GRS} and Radziszowski    \cite{survey} for  a dynamic survey.\medskip

In this paper we investigate  Ramsey numbers of odd cycles and odd cycles  with chords under Gallai colorings,   
 where a \dfn{Gallai coloring} is a coloring of the edges of a complete graph without rainbow triangles (that is, a triangle with all its edges colored differently). Gallai colorings naturally arise in several areas including: information theory~\cite{KG}; the study of partially ordered sets, as in Gallai's original paper~\cite{Gallai} (his result   was restated in \cite{Gy} in the terminology of graphs); and the study of perfect graphs~\cite{CEL}. There are now a variety of papers  which consider Ramsey-type problems in Gallai colorings (see, e.g., \cite{C9C11, C7, chen, c5c6,GS, exponential, Hall,  K4, C6C8,W4}).    More information on this topic  can be found in~\cite{FGP, FMO}.  \medskip
 
A \dfn{Gallai $k$-coloring} is a Gallai coloring that uses at most $k$ colors. 
 Given an integer $k \ge 1$ and graphs $H_1,  \ldots, H_k$,  \dfn{Gallai-Ramsey number}  $GR(H_1,  \ldots, H_k)$ is defined to be  the least integer $N$ such that every Gallai $k$-coloring of $K_N$   contains a monochromatic copy of $H_i$ in color $i$ for some $i \in [k]$. When $H = H_1 = \dots = H_k$, we simply write $GR_k(H)$.    Clearly, $GR_k(H) \leq R_k(H)$ for all $k\ge1$ and $GR(H_1, H_2) = R(H_1, H_2)$.     
Theorem~\ref{general} below is a result of Gy\'{a}rf\'{a}s,   S\'{a}rk\"{o}zy,  Seb\H{o} and   Selkow~\cite{exponential} which characterizes    the general behavior of $GR_k(H)$.

\begin{thm} [\cite{exponential}]\label{general}
Let $H$ be a fixed graph  with no isolated vertices 
 and let $k\ge1$ be an integer. Then
$GR_k (H) $ is exponential in $k$ if  $H$ is not bipartite,    linear in $k$ if $H$ is bipartite but  not a star, and constant (does not depend on $k$) when $H$ is a star.				
\end{thm}

It turns out that for some graphs $H$ (e.g., when $H=K_3$),  $GR_k(H)$ behaves nicely, while the order of magnitude  of $R_k(H)$ seems hopelessly difficult to determine.  We will utilize the following important structural result of Gallai~\cite{Gallai}.

\begin{thm}[\cite{Gallai}]\label{Gallai}
	Let  $(G, \tau) $  be a Gallai $k$-colored complete graph    with   $|V(G)|\ge  2 $.  Then  $V(G)$ can be partitioned into nonempty sets  $V_1,  \dots, V_p$ with $p\ge2$ so that    at most two colors are used on the edges in $E(G)\less (E(V_1)\cup \cdots\cup  E(V_p))$ and only one color is used on the edges between any fixed pair $(V_i, V_j)$ under $\tau$, where $E(V_i)$ denotes the set of edges with both ends in  $ V_i $ for all $i\in [p]$. 
\end{thm}

 The partition $\{V_1,   \dots, V_p\}$ given in Theorem~\ref{Gallai} is  a \dfn{Gallai partition} of    $(G, \tau)$.  Let $(\mathcal{R},\tau)$ be obtained from  $(G, \tau)$ by first contracting each $V_i$ into a single vertex $v_i$ and then coloring $v_iv_j$ by the unique color on the edges between $V_i$ and $V_j$ in $(G, \tau)$.   We say  $(\mathcal{R},\tau)$ is   the  \dfn{reduced graph} of $(G,\tau)$ corresponding to the  Gallai partition $\{V_1,   \dots, V_p\}$. Note that   $\mathcal{R}=K_p$.  
By Theorem~\ref{Gallai},  all edges in $\mathcal{R}$ are colored by at most two colors under $\tau$.  One can see that any monochromatic copy of $H$ in $(\mathcal{R},\tau)$  will result in a monochromatic copy of $H$ in $(G, \tau)$. It is not  surprising  that    the  $2$-color  Ramsey number  $R_2(H)$ plays an important role in determining the value of  $GR_k( H)$ when $H$ is a complete graph.    Fox,  Grinshpun and  Pach~\cite{FGP} posed the following  conjecture.  

\begin{conj}[\cite{FGP}]\label{Fox} For all  $k\ge1$ and $t\ge3$,
\[
GR_k( K_t) = \begin{cases}
			(R_2(K_t)-1)^{k/2} + 1 & \text{if } k \text{ is even} \\
			(t-1)  (R_2(K_t)-1)^{(k-1)/2} + 1 & \text{if } k \text{ is odd.}
			\end{cases}
\]
\end{conj}

  The first case of Conjecture \ref{Fox} follows directly from a result of Chung and Graham~\cite{chgr} in 1983. A simpler proof of the first case of Conjecture \ref{Fox} can be found in~\cite{exponential}. The next open case, when $t=4$,  was recently settled in~\cite{K4}.   \medskip

In this paper, we    focus on determining the exact values of  $GR_k(\Theta_{2n+1})$  for all $k\ge1$ and  $n\ge3$,  and $GR_k(C_{2n+1})$ for all $k\ge1$ and $n\ge2$.   It is worth noting that the construction of Bondy and Erd\H{o}s mentioned earlier for  $R_k(C_{2n+1})$ contains no rainbow triangles. Hence  $GR_k(\Theta_{2n+1})\ge n\cdot 2^k+1$ and $GR_k(C_{2n+1})\ge n\cdot 2^k+1$ for all $k\ge1$ and $n\ge2$.  More recently,   the exact value of $GR_k(C_{2n+1}) $ for $2\le n\le7$    has been completely settled, see  \cite{ C9C11, C13C15, C7, c5c6};   the exact value of $GR_k(\Theta_{5}) $   was recently settled in \cite{Theta5}.   We determine the exact value of $GR_k(\Theta_{2n+1}) $ for all $k\ge 1$ and $n\ge3$.  We include $C_5$ in Theorem~\ref{oddchorded} below  in order to   provide     a unified proof for the Gallai-Ramsey number  of all odd cycles on at least five vertices. 
 
\begin{thm}\label{oddchorded}
For   all     $k\geq 1$,   $GR_k(C_{5})=2\cdot 2^k+1$ and $GR_k(\Theta_{2n+1})=n\cdot 2^k+1$  for all  $n\ge3$. 

\end{thm}

Note that  every monochromatic copy of $\Theta_{2n+1}$ contains a monochromatic copy of $C_{2n+1}$. Theorem~\ref{oddchorded} not only   implies  the exact value for $GR_k(C_{2n+1})$  but also provides   a unified proof for the Gallai-Ramsey number  of $C_{2n+1}$ for all $k\ge1$ and $n\ge2$. 

\begin{thm}\label{odd}
For   all     $k\geq 1$ and   $n\geq 2$,  $GR_k(C_{2n+1})=n\cdot 2^k+1.$

\end{thm}

 Theorem~\ref{odd} implies that  Conjecture~\ref{OddCycle} is true under Gallai colorings.  
 We want to point out here that Wang et al. have also posted a paper \cite{falseclaim}  claiming the   result of Theorem~\ref{odd}.  We were thus compelled to share our independent proof  in the manuscript \cite[Theorem 1.7]{earlyversion} which was not submitted for publication. 
 Our    proof of Theorem~\ref{oddchorded} given in Section~\ref{OD}   relies   on    an  upper bound result for   Gallai-Ramsey number  of even cycles   in \cite[Theorem 1.6]{evencycle}. However, the proof of \cite[Theorem 1.6]{evencycle} is rather long. In this paper we prove Theorem \ref{oddchorded}   using a more relaxed upper bound on $GR_k(C_{2n})$ given in Theorem \ref{even} below. We include a short proof of Theorem \ref{even} here for completeness. We need to introduce more notation before we state the result.\medskip
 
 For the remainder of the paper, we use $(G, \tau)$ to denote a Gallai $k$-colored complete graph, where $G$ is a complete graph and   $\tau: E(G)\rightarrow [k]$ is a Gallai $k$-coloring.   For  each $(G, \tau)$, let  $G^\tau_i$   denote the spanning subgraph of $G$ with   $E(G^\tau_i): = \{e\in E(G)\mid  \tau(e)=i\}$ for all $i\in[k]$. We simply write $G^\tau_r$ if the color $i$ is red; $G^\tau_b$ if the color $i$ is blue.   
For  every   $n\in \mathbb{N}$,    let  $q_\tau(G, n)$ denote   the number of colors $i\in[k]$ such that  $G^\tau_i$ has  a component of order at
least $n$.  Then  $q_\tau(G, n)\le k$.   We begin with  Lemma~\ref{q_tau} (which follows directly from the proof of Lemma 9 given in \cite{Hall}).   A proof of Lemma\ref{q_tau} can be found in \cite[Lemma 1.5]{evencycle}.  

\begin{lem}[\cite{Hall}]\label{q_tau}
Let  $(G, \tau) $  be a Gallai $k$-colored complete graph    with   $|V(G)|\ge n\ge2 $.  
 Let $\{V_1, \ldots, V_p\}$ be  a Gallai partition of $(G, \tau)$   with  $p \ge 2$ as small as possible. Then   $q_\tau(G, n)\ge1$, and  for any color $i$ on the edges  in $(\mathcal{R}, \tau)$, $\mathcal{R}^\tau_i$   is connected.   
  
\end{lem}

 Theorem \ref{even} below establishes an upper bound for  $GR_k(C_{2n})$ for all $k\ge1$ and  $n\ge2$.

\begin{thm}\label{even}
 Let  $(G, \tau) $  be a Gallai $k$-colored complete graph    with   $|V(G)|\ge n\ge 2  $. If  
\[|V(G)| \ge  (n-1)\cdot q_\tau(G,  n)+2n+2,\]  then $(G, \tau)$ has a   monochromatic       copy of      $ C_{2n}$. 
\end{thm}

We conclude this section by introducing more notation and list several known results that will be used in the proof  of Theorem~\ref{even} and Theorem~\ref{oddchorded}.     Given a graph $G$ and a  set  $S\subseteq V(G)$,  we use   $|G|$    to denote  the  number
of vertices    of $G$,      $G\less S$ the subgraph    obtained from $G$ by deleting all vertices in $S$,  and $G[S] $    the  subgraph    obtained from $G$ by deleting all vertices in $V(G)\less S$. 
      For  two disjoint sets $A, B\subseteq V(G)$,    $A$ is \dfn{complete} to $B$ in $G$  if each vertex in $A$ is adjacent to all vertices in  $B$, and \dfn{anti-complete} to $B$ in $G$  if no  vertex in $A$ is adjacent to any vertex  in  $B$.   
 Let $(G,\tau)$ be a Gallai $k$-colored complete graph.   
For  two disjoint sets $A, B\subseteq V(G)$,   $A$ is \dfn{mc-complete} to $B$   if all the edges between $A$ and $B$  in    $(G,\tau)$ are colored the same color. We simply say      $A$ is     \dfn{$j$-complete} to $B$    if all the edges between $A$ and $B$  in   $(G,\tau)$ are colored by some color $j\in[k]$,   and   $A$ is \dfn{blue-complete}     to $B$    if all the edges between $A$ and $B$  in $(G,\tau)$ are colored  blue.  We say a vertex $x\in V(G)$ is \dfn{blue-adjacent} to a vertex $y\in V(G)$ if the edge $xy$ is colored blue in $(G,\tau)$, and $x$ is \dfn{blue-complete} to an edge $yz\in E(G)$ if $x$ is  blue-complete to $\{y, z\}$ in $(G,\tau)$.  
Similar definitions hold when blue is replaced by another color. For convenience, we use  $A \less B$ to denote  $A-B$; and  $A \less b$ to denote  $A \less \{b\}$ when $B=\{b\}$. 
We use the convention   ``$S:=$'' to mean that $S$ is defined to be the right-hand side of the relation.

 \begin{thm}[\cite{BE, theta}]\label{2n+1}     For all  $n\geq 2$, $R_2(C_{2n+1})=R_2(\Theta_{2n+1})=4n+1$. 
\end{thm}

 For a bipartite graph $G$ with bipartition $(M, N)$, let     $\delta(N):= \min\{d_G(x): x\in N\}$ and $\Delta(N):= \max\{d_G(x): x\in M\}$..  Lemma~\ref{Hal2}   follows from Lemma 12, Lemma 13, Lemma 14 in \cite{Hall}. 
 
\begin{lem}[\cite{Hall}]\label{Hal2}
Let $G$ be a bipartite graph  with bipartition $(M, N)$ such that $|N| \geq 2$, $|M| \geq 4$ and  $\delta(N) \ge (|M| + 1)/{2}$. Then $G$ has a cycle of length $ 2\ell $ for any $\ell$ satisfying $2 \leq \ell \leq \min\{|N|, \delta(N)-1\}$,  or $\Delta(N)=(|M|+1)/2$, $ M=M_1\cup  M_2\cup M_3$ and $N=N_1\cup N_2$, where $M_1, M_2, M_3, N_1,  N_2 $ are non-empty, pairwise disjoint  sets in $G$, $|M_3|=1$, $|M_1|=|M_2|$, and $N_i  $ is  complete to $M_i\cup M_3$ but anti-complete to $M_{3-i}$ for all $i\in[2]$.    \end{lem}

\begin{lem}[\cite{C9C11}]\label{L2}

 Let $\{V_1,\ldots,V_p\}$ be a Gallai partition of $(G, \tau)$ with $|V_1| \le \cdots\le |V_p|$ and $|G| \ge 2n + 1$. If $|V_p| \le n$ and the corresponding reduced
graph $\mathcal{R}$ is monochromatic, say
blue, then $(G, \tau )$ contains a blue $C_{2n+1}$.
\end{lem}

\begin{lem}[\cite{C9C11}]\label{L}
Let    $Y, Z$ be   two disjoint sets of vertices in $(G, \tau)$ with $|Y|  \ge   n$ and $ |Z| \ge n$. If $(G, \tau)$ is  $C_{2n+1}$-free and   $Y$ is mc-complete, say blue-complete, to  $Z$,  then no vertex  in $ V(G) \setminus (Y \cup Z)$ is blue-complete to $Y \cup Z$ in $(G, \tau)$.  Moreover, if $|Z| \ge n+1$, then $(G[Z],\tau)$   has no blue edges. Similarly, if  $|Y| \ge n+1$, then   $(G[Y],\tau)$ has no blue edges.
\end{lem}

\section{ Proof of Theorem  \ref{even}}\label{EC}

   Let $(G,\tau)$ and $k$ be as given in the statement.  Suppose $(G,\tau)$ is $C_{2n}$-free. Choose $(G,\tau)$ so that $q_\tau(G, n)$  is minimum. 
      Let $q:=q_\tau(G,n)$. By Lemma~\ref{q_tau}, $q\ge1$. 
   We may assume that for each  color $i\in[q]$,  $G_i^\tau$ has a component of order at
least $n$.  
Let $X_1,   \ldots, X_q$  be   disjoint subsets of $V(G)$ such that   for each  $i \in [q]$,   $X_i$  (possibly empty) is  mc-complete in color $i$ to $V(G) \less \bigcup_{i=1}^q X_i$.  Choose $X_1, \ldots, X_q$ so that $|G|-\sum_{i=1}^q |X_i| \ge n$ and $\sum_{i=1}^q |X_i| $ is as large as possible.  Denote $X:= \bigcup_{i=1}^q X_i$.  Then $|G\less X|\ge n\ge2$.  Since $(G,\tau)$ has no rainbow triangle, we see that for all $i, j\in[q]$ with $i\ne j$, all the edges between $X_i$ and $X_j$ are colored by color $i$ or color $j$. We next  prove  a series of  claims.

\begin{claim}\label{evenX_i} 
For all $i \in [q]$, $|X_i| \le n-1$.
\end{claim}
\pf Suppose $|X_i|\ge n$ for some color $i\in[q]$, say blue. By the choice of $X_1, \ldots, X_q$, $X_i$ is blue-complete to $V(G)\less X$. Thus  $(G,\tau)$ has a blue    $C_{2n}$ using edges between $X_i$ and $V(G)\less X$, a contradiction. \qed\medskip

Let  $V_1, \ldots, V_p$ be a Gallai partition of $(G \less X,\tau)$ with $p\ge 2$ as small as possible.  We may assume that $|V_1|  \le \cdots \le |V_p|$.  By Theorem~\ref{Gallai} and Lemma~\ref{q_tau},  all the edges of the reduced graph of 
$(G \less X,\tau)$ are colored by at most two colors in $[q]$, say red or blue.  Then for all $i\in [p-1]$, $V_i$ is either  red- or  blue-complete to $V_p$ in $(G\less X, \tau)$.  Let 
\[\begin{split}
\mathcal{V}_r &:=    \{V_i  \in \{V_1, \ldots, V_{p-1}\} \mid V_i \text{ is   red-complete to } V_p \text{ in } (G\less X, \tau)\} \text{ and }\\
\mathcal{V}_b &:=    \{V_i  \in \{V_1, \ldots, V_{p-1}\} \mid V_i \text{ is   blue-complete to } V_p \text{ in } (G\less X, \tau)\}.
\end{split}\]
Let $R=\bigcup_{V_j\in \mathcal{V}_r}V_j$ and $B=\bigcup_{V_j\in \mathcal{V}_b}V_j$. Then $R\cup B=V(G)\less (X\cup V_p)$,  and $R$ and $B$ are disjoint. We may further assume that $X_1$ is red-complete to
$V(G) \less X$ and $X_2$ is blue-complete to $V(G) \less X$.

\begin{claim}\label{evenB} 
$\min\{|B|, |R|\}\ge 1$ and so $q\ge2$.
\end{claim} 

\pf Suppose $\min\{|B|, |R|\}= 0$. We may assume that $B=\emptyset$. Since $p\ge 2$, we see that $|R|\ge1$ and $R$ is red-complete to $V_p$. Then $|V_p|\le n-1$, else, let  $X'_1:=X_1\cup R$ and $X'_i:=X_i$ for all $i\in\{2, \ldots, q\}$.  But then  
$ |X'_1|+\cdots+|X_q'| =|X \cup R | >|X|$,
 contrary to  the  choice of $X_1, \ldots, X_q$.  Similarly, $|R|\le n-1$. By Claim \ref{evenX_i}, $|X_i|\le n-1$ for all $i\in [q]$. But then $|G|=|X|+|V_p|+|R|\le (n-1)q+2(n-1)$, contrary to the assumption that $|G| \ge  (n -1)q + 2n+2$.  This proves that $\min\{|B|, |R|\}\ge 1$. By the choice of   $p$ and Lemma~\ref{q_tau}, both $G_r^\tau\less X$ and  $G_b^\tau\less X$ are connected. Thus $q\ge2$. \qed

\begin{claim}\label{evenVpp} 
$|V_p|\le n-1$ and so $X_i\ne \emptyset $ for   any  color $i\in[q]$ that is neither red nor blue.
\end{claim} 

\pf Suppose $|V_p|\ge n$.  Then every vertex in $B \cup R $ is either  red- or  blue-complete to $V_p$.  Let $X'_1:=X_1\cup R$, $X'_2:=X_2\cup B $, and $X'_i:=X_i$ for all $i\in\{3, \ldots, q\}$.  But then  
$ |X'_1|+\cdots+|X_q'| =|X \cup B  \cup R | >|X|$,
 contrary to  the  choice of $X_1, \ldots, X_q$.   This proves that $|V_p|\le n-1$. Next,
suppose there exists a color $j\in[q]$ such that $j$ is neither red nor blue but $X_j = \emptyset$. Then no edges
between pairs of $X_1, \ldots,X_q$ are colored by color $j$ in $(G,\tau)$. But then $G_j^\tau$
has no component of
order at least $n$, because $|V_\ell |\le |V_p|\le n-1$ for all $\ell\in[p-1]$, and $|X_i
| \le n-1$ for all $i \in [q]$ by Claim \ref{evenX_i}, a contradiction.   \qed
  
\begin{claim}\label{evenq} 
 $n\ge 3$ and 
$q=2$.
\end{claim} 

\pf  Suppose   $n=2$.  By Claim \ref{evenX_i} and   Claim \ref{evenVpp}, $|V_p|=1$ and $|X_i|\le1$ for all $i\in [q]$.  Then $|G\less \bigcup_{i=3}^q X_i|\ge q+6- (q-2)=8$. Thus   $(G\less \bigcup_{i=3}^q X_i,\tau)$ contains a red or blue $C_4$ because $R_2(C_4)=6$ \cite{CS}, a contradiction. 
Suppose   next $n\ge 3$ and 
$q\ne 2$. By Claim~\ref{evenB}, $q\ge 3$ and the color $q$ is neither red nor blue. By Claim~\ref{evenVpp},   $X_q\ne \emptyset$ and $|V_p|\le n-1$.  Thus $q_\tau(G\less X_q,n)=q_\tau(G,n)-1=q-1$.   By Claim~\ref{evenX_i}, $|X_q|\le n-1$, and so $|G\less X_q|\ge (n-1)q+2n+2-(n-1)=(n-1)\cdot q_\tau(G\less X_q,n)+2n+2$. By the minimality  of $q$, $(G\less X_q, \tau)$ has a monochromatic copy of  $C_{2n}$, a contradiction.  
\qed\medskip

By Claim~\ref{evenq}, $q=2$  and $n\ge 3$ 
and  so $|G|\ge 2(n-1)+2n+2=4n$. We may further assume that $|X_1\cup R|\ge |X_2\cup B|$. 
Then $|X_1\cup R|\ge \left\lceil(|G|-|V_p|)/2\right\rceil\ge n+2$ and    
$|V_p\cup X_1\cup R| \ge |V_p|+ \left\lceil(|G|-|V_p|)/2\right \rceil\ge 2n+1$.  Let 
\[
\begin{split}
W:=  
 \begin{cases}  
X_1\cup R  &  \text{if } |X_1\cup R|\le 2n+1 \\
X_1\cup  \left(\bigcup_{V\in \mathcal{V}_r'€™}V\right)      &  \text{if } |X_1\cup R|> 2n+1,\\
\end{cases} 
\end{split}
\]
where 
 $ \mathcal{V}_r' \subseteq  \mathcal{V}_r$  is chosen  so that  $ |W| \le 2n+1$ is  as large as possible when $|X_1\cup R|> 2n+1$.    In both cases, we have $|V_p\cup W|\ge 2n+1$ because in the second case $|V_p|\ge |V_i|$ for all $i\in[p-1]$.   Let $V_p^*\subseteq V_p$ be such that $|V^*_p\cup W|= 2n+1$. Let $M:=V^*_p\cup W$ and $N:= V(G)\less M$. Then  $|N|\ge 4n-(2n+1)=2n-1$ and $N=(V_p\less V^*_p)\cup B\cup (R\less W)\cup X_2$. By   Claim~\ref{evenVpp}, $|W|=|M\less V_p|\ge (2n+1)-(n-1)=n+2$.    It follows that   every vertex in $V_p\less V^*_p$ is red-adjacent to    every vertex in $W$, and  every vertex in $N\less V_p$ is either red- or blue-adjacent to  at least  $\lceil|M|/2\rceil= n+1$ vertices in $M$.  Let 
\[
\begin{split}
N_r &:=    \{v  \in N \mid v \text{ is   red-adjacent to at least } n+1 \text{ vertices in } M \}  \text{ and }\\
N_b &:=    \{v  \in N \mid v \text{ is   blue-adjacent to at least } n+1 \text{ vertices in }  M\}.
\end{split}
\]
Then $|N_r|+ |N_b|=|N|\ge 2n-1$.  Since every vertex in $V_p\less V^*_p$ is red-adjacent to  at least  $ |W|\ge n+2$ vertices in $M$, we see that   $N_b\cap V_p=\emptyset$ and $V_p\less V^*_p\subseteq N_r$.   For each color $j\in \{r, b\}$,    let $H_j$   be the bipartite subgraph of $G$ with   bipartition  $(M, N_j)$      such that  $E(H_j)$ consists  of all     edges in color $j$ between    $M$ and $N_j$  in $(G,\tau)$.
Then $H_j$ is $C_{2n}$-free because $(G,\tau)$ is $C_{2n}$-free. \medskip

Suppose  $|N_b|\ge n$. Since $N_b\cap V_p=\emptyset$, we see that every edge  between $M$ and $N_b$ is  colored red or blue in $(G, \tau)$.  Note that  $\min\{d_{H_b}(v): v\in N_b\}\ge n+1=(|M|+1)/2$ and  $H_b$ is $C_{2n}$-free. By Lemma \ref{Hal2} applied to $H_b$ with $\ell:=n$, we see that  $ M=M_1\cup  M_2\cup M_3$ and $N_b=N_1\cup N_2$, where $M_1, M_2, M_3,N_1,  N_2$ are non-empty, pairwise disjoint  sets in $H_b$, $|M_3|=1$, $|M_1|=|M_2|=n$, $N_i  $ is  blue-complete to $M_i\cup M_3$ and  red-complete to $M_{3-i}$ for all $i\in[2]$.
 Since    $G$ contains no rainbow triangle, every edge between $M_1$ and $M_2$  is  colored   red or blue.  
Let  $x_1, y_1,   z_1 \in M_1$  and $x_2, y_2, z_2\in M_2$ be all distinct. Since each edge between $M_1$ and $M_2$  is  colored   red or blue, we may further assume that $x_1x_2$ and $y_1y_2$ are colored the same. For each  $i\in[2]$, let $Q_i^r$ be an $(x_i, y_i)$-path between  $M_i$ and $N_{3-i}$  such that $|Q_1^r|+|Q_2^r|=2n$; $Q_i^b$ be an $(x_i, y_i)$-path between  $M_{i}$ and $N_{i}$ such that $|Q_1^b|+|Q_2^b|=2n$.   But then $(G,\tau)$ contains  a  red $C_{2n}$ with edge set $E(Q_1^r)\cup E(Q_2^r)\cup\{x_1y_1, x_2y_2\}$ when $x_1x_2$ and $y_1y_2$ are colored red, or   a  blue $C_{2n}$ with edge set $E(Q_1^b)\cup E(Q_2^b)\cup\{x_1y_1, x_2y_2\}$ when $x_1x_2$ and $y_1y_2$ are colored blue,   a contradiction.  \medskip
 
It remains to consider the case  $|N_b|\le n-1$. Then    $|N_r|\ge n$ and  $\min\{d_{H_r}(v): v\in N_r\}\ge n+1=(|M|+1)/2$.  Note that $H_r$ is $C_{2n}$-free. By Lemma \ref{Hal2} applied to $H_r$ with $\ell=n$, we see that  $ M=M_1\cup  M_2\cup M_3$ and $N_r=N_1\cup N_2$, where $M_1, M_2, M_3,N_1,  N_2$ are non-empty, pairwise disjoint  sets in $H_r$, $|M_3|=1$, $|M_1|=|M_2|=n$, $N_i  $ is  red-complete to $M_i\cup M_3$ and  no edges between $N_i$ and   $M_{3-i}$  are colored red for each  $i\in[2]$.  Since every vertex in $V_p\less V^*_p$ is red-adjacent to  at least  $  n+2$ vertices in $M$, it follows that    $N_r\cap V_p=\emptyset$.  By the choice of $W$,    $N_i$ is blue-complete to $M_{3-i}$ for each  $i\in[2]$. 
 Note that  every edge between $M_1$ and $M_2$  is  colored   red or blue  because      $(G,\tau)$ contains no rainbow triangle.    Similar  to the  argument   in the previous  paragraph, we obtain a monochromatic $C_{2n}$ in $(G,\tau)$,  a contradiction. \medskip
 
   This completes the proof of Theorem~\ref{even}. \qed

\section{Proof of Theorem \ref{oddchorded}}\label{OD}

Let $n$, $k$ be as given in the statement. It suffices to show that $GR_k(C_5)\le 2\cdot 2^k+1$ for all $k\ge1$ and $GR_k(\Theta_{2n+1})\le n\cdot 2^k+1$ for all $k\ge1$ and $n\ge3$. The case when $k=1$ is trivial and the case when $k=2$ follows from Theorem \ref{2n+1}. So we may assume that $k\geq 3$. Let $G:=K_{n\cdot 2^k+1}$ and let $\tau:E(G)\rightarrow [k]$ be any  Gallai $k$-coloring of $G$. 
Suppose that  $(G,\tau)$ is $\Theta_{2n+1}$-free when $n\ge3$ and $C_5$-free when $n=2$. Choose $(G,\tau)$ with $k$ minimum.  Then $\tau$ is an onto mapping. 	 
Let $X_1,   \ldots, X_k$  be   disjoint subsets of $V(G)$ such that   for each  $i \in [k]$,   $X_i$  (possibly empty) is  mc-complete in color $i$ to $V(G) \less \bigcup_{i=1}^k X_i$ under $\tau$.  Choose $X_1, \ldots, X_k$ so that $\sum_{i=1}^k |X_i| \le (k+1)n$ is as large as possible.  Let $X:= \bigcup_{i=1}^k X_i$.  Then $|X|\le (k+1)n$ and $|G\less X|\ge n\cdot 2^k+1-(k+1)n\ge 4n+1$. Since $(G,\tau)$ has no rainbow triangle, we see that for all $i, j\in[k]$ with $i\ne j$, all the edges between $X_i$ and $X_j$ are colored by color $i$ or color $j$. We next  prove several  claims.

\begin{claim}\label{noC2n+1} $(G,\tau)$ is $C_{2n+1}$-free for all $n\ge2$. 
\end{claim}

\pf  Suppose $(G,\tau)$ contains a monochromatic    $C:=C_{2n+1}$,  
 say    with vertices $x_1,x_2,\ldots,x_{2n+1}$ in order. Then $n\ge3$. We may assume that all edges of $C$ are colored blue.  Then  no chord of $C$ is colored blue   because  $(G,\tau)$ is  $\Theta_{2n+1}$-free. We may further assume that $x_1x_3$ is colored red. Then $x_1x_j$ is colored red for all $j\in \{3, 4, \ldots, 2n\}$ because $(G, \tau)$ has a no rainbow triangle. It follows that all chords of  $C$  are colored red.  Let $H$ be the graph with $V(H)=V(C)$ and $E(H)$ consisting of all chords of $C$. Then $H$ is the complement of $C_{2n+1}$. It can be easily checked that   $H$  contains a chorded $C_{2n+1}$ because $n\ge3$.   Thus $(G,\tau)$ contains a red copy of   $\Theta_{2n+1}$, a contradiction. \qed\medskip
 
\begin{claim}\label{Odd_X_i bound}
For all $i \in [k]$, $|X_i| \le n-1$.
\end{claim}
\pf Suppose $|X_j| \ge n  $ for some color $j \in [k]$. Since $|G\less X|\ge 4n+1$, by Lemma \ref{L} applied to $X_j$ and $V(G)\less X$,  we see that $(G\less X,\tau)$ has no  edge in color $j$. By minimality of $k$, $|G\less X|\le n\cdot 2^{k-1}$. But then     
\[
|G| = |X|+|G \less X| \le (k+1)n+n\cdot 2^{k-1} <n\cdot 2^k+1
\] 
for all $k\ge3$ and $n\ge 2$, a contradiction. \qed

\begin{claim}\label{Odd_empty} 
 $X_i = \es$ for   some  $i \in [k]$, and so $|X| \le (k-1)(n-1)$.
 \end{claim}
\pf   Suppose $X_i \ne \es$ for every $i \in [k]$.  By Claim~\ref{Odd_X_i bound}, $|X| \le k(n-1)$.  Then 
\[
|G \less X| \ge n\cdot 2^k+1 -  k(n-1)  > (n-1)k + 2n+2,
\]
for all $k \ge 3$ and $n \ge 2$. By Theorem~\ref{even} and the fact that $q_\tau(G\less X, n)\le k$, we see that $(G\less X,\tau)$     contains  a monochromatic  copy of    $C_{2n}$, say in color 1. Since $X_1\neq \emptyset$, we see that   $(G,\tau)$ contains  a monochromatic copy of   $C_{2n+1}$  in color 1,  contrary to Claim~\ref{noC2n+1}. This proves that $X_i = \es$ for   some  $i \in [k]$. By  Claim  \ref{Odd_X_i bound},  $|X| \le (k-1)(n-1)$. \qed\medskip

\begin{figure}[h]
\centering
\includegraphics[scale=0.5]{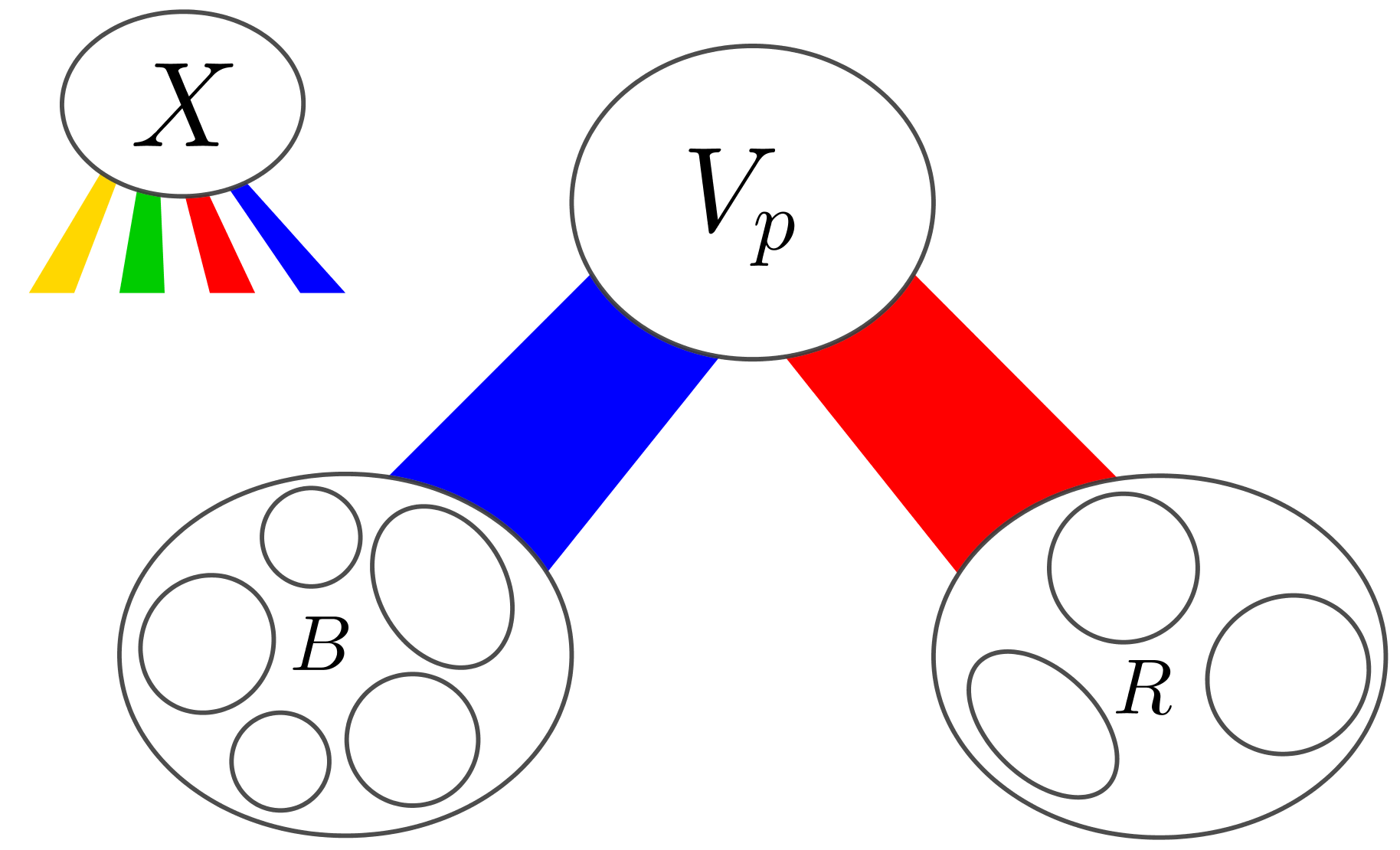}
\caption{An overview of $(G, \tau)$.}
\label{f1}
\end{figure}

Let $\{V_1, \ldots, V_p\}$ be a Gallai partition of $(G \less X,\tau)$  with  $p \ge 2$ and $|V_1|  \le \cdots \le |V_p|$.  By Theorem~\ref{Gallai}, we may assume that  every  edge  of the corresponding reduced graph $(\mathcal{R},\tau)$ of  $(G \less X,\tau)$ is  colored either red or blue.   By the choice of $X_1, \ldots, X_k$, we may further assume that  $X_1$ is red-complete to $V(G)\less X$ and $X_2$ is blue-complete to $V(G)\less X$.     By Theorem~\ref{2n+1},    $p \le 4n$.  Thus $|V_p|\ge2$.   Note that  for all $i\in [p-1]$, $V_i$ is either  red- or  blue-complete to $V_p$ in $(G, \tau)$.  Let 
\[
\begin{split}
 \mathcal{V}_r   &:=    \{V_i  \in \{V_1, \ldots, V_{p-1}\} \mid V_i \text{ is   red-complete to } V_p \text{ in } (G,\tau)\}\text { and }\\
 \mathcal{V}_b  &:=    \{V_j  \in \{V_1, \ldots, V_{p-1}\} \mid V_j \text{ is   blue-complete to } V_p \text{ in } (G,\tau)\}.
\end{split}
\]
Let $R:=\bigcup_{V_j\in  \mathcal{V}_r  }V_j$ and $B:=\bigcup_{V_j\in \mathcal{V}_b }V_j$.  Then $R\cup B=V(G)\less (X\cup V_p)$.  
An overview of $(G,\tau)$ is
depicted in Figure \ref{f1}.

\begin{claim}\label{Odd_B_R}
$|R \cup B| \ge 2n + 2$.
\end{claim}
\pf Suppose  $|R \cup B| \le 2n+1$. Since $p\ge 2$, we see that $|R \cup B|\ge 1$.   
 Let $X'_1:=X_1\cup R$, $X'_2:=X_2\cup B$, and $X'_i:=X_i$ for all $i\in\{3, \ldots, k\}$. Then $V_p=V(G)\less\bigcup_{i=1}^k X_i'$, and for all $i\in[k]$, $X_i'$ is mc-complete to $V_p$ in $(G,\tau)$. Thus
\[
\sum_{i=1}^k |X_i|<\sum_{i=1}^k |X'_i|=|X \cup R \cup B| \le (k-1)(n-1)  + 2n+1 = (k+1) n - (k-2) <   (k+1)n,
\]
contrary to  the  choice of $X_1, \ldots, X_k$.   \qed

\begin{claim}\label{Odd_Vp}
$|V_p| \ge n+1$.
 \end{claim}
\pf Suppose $|V_p| \le n$.  We may assume that $|R|\ge |B|$.   Then 
\[|R|+ |B|= |G|-|V_p|-|X| \geq  (n\cdot 2^k+1)-n-(k-1)(n-1)=n\cdot 2^k-kn +k\ge 5n+3.\] 
 Thus $|R|>2n+1.$  By Lemma \ref{L2} and Claim~\ref{noC2n+1}, $(G[R],\tau)$ must have both red and blue edges. Suppose $|X_1\cup V_p|\ge n$. By Lemma \ref{L} applied to $(G,\tau)$ with $Y=R$ and $Z=X_1\cup V_p$, we see that $(G[R],\tau)$ has no red edges, a contradiction. Thus $|X_1\cup V_p|\le n-1$ and so $|V_p|\le n-1$. It follows that
for any color $i\in [k]$ other than red and blue, $G_i^\tau\less X$ has no component of order at least $n$. Thus
$q_\tau (G\less X, n) \le 2$ and then $|G\less X| > (n-1)\cdot q_\tau (G\less X, n) +2n+2$. By Theorem \ref{even}, $(G\less X, \tau)$ contains a red or blue $ C_{2n}$. This implies that either $X_1=\emptyset$ or $X_2=\emptyset$.
\vskip 2mm

Suppose   $|R|\ge 4n $.  Let $  \mathcal{V}_r'  \subseteq   \mathcal{V}_r $  be chosen  so that  $   |\bigcup_{V_i\in  \mathcal{V}_r' }V_i|\ge4n$ is  as small as possible. Let  $U:=    \bigcup_{V_i\in  \mathcal{V}_r'  }V_i $.    Then  $|U|\leq 5n-2$ because $|V_1|\le \cdots \le|V_p| \le n-1$.    Since $q_\tau(G\less X,n)\le2$, we see that $q_\tau(G[U],n)\le2$. Then $|U|\ge 4n\ge (n-1)\cdot q_\tau(G[U],n)+2n+2$. By Theorem \ref{even},  $(G[U], \tau)$ contains a monochromatic  copy of  cycle $C$ of length $2n$, say with vertices $u_1, u_2, \ldots, u_{2n}$ in order. Let $U_1=\{u_1,u_3,\ldots,u_{2n-1}\}$ and  $U_2=\{u_2,u_4,\ldots,u_{2n}\}$. Since $U$ is red-complete to $V_p$, we see that the cycle $C$ must be blue. Then $X_2=\emptyset$. Let $x,y\in V_p$ be distinct.   Then 
\[
 |V(G)\less (X\cup U\cup \{x,y\})|\ge n\cdot 2^k+1-(k-1)(n-1)-(5n-2)-2\ge n+3.
 \]
Let  $W\subseteq V(G)\less (X\cup U\cup \{x,y\})$ with $|W|=n$.  By the choice of $U$, every edge between $W$ and $V(C)$  is colored   either red or blue in $(G,\tau)$. Note that  no  vertex in $ W$ is   blue-adjacent to   any edge on the cycle $C$, else, we obtain a blue   $C_{2n+1}$, contrary to Claim~\ref{noC2n+1}.  It follows that each vertex in $ W$ is red-adjacent to at least $n$ vertices on $C$. 
Let 
\[
\begin{split}
W_0: &=\{v\in W  \mid v \text{ is red-adjacent to at least  } n+1 \text{  vertices on  }  C\},\\
W_1: &=\{v\in W\less W_0  \mid v \text{ is red-complete to } U_1  \text{  only in   } G\} \text{ and } \\
W_2: &=\{v\in W\less W_0 \mid v \text{ is red-complete to } U_2  \text{  only in   } G\}.\\
\end{split}
\]
Then  $W=W_0\cup W_1\cup W_2$ and   $W_i$ is blue-complete to $U_{3-i}$ for each $i\in[2]$.   Let $H$ be the
bipartite subgraph of $G$ with bipartition  $(V(C), W)$ such that $E(H)$ consists of all the red edges between $W$ and $V(C)$ in $(G,\tau)$.  Suppose $H$ has  a  path $P$ on $2n-3$ vertices     with both ends in  $V(C)$.  We may assume that $u_1, u_j$ are the ends of $P$ for some $j\in[2n]$ with $j\ne 1$.  
Since  $|U\cap V(P)|=n-1$, we see that $|U\less V(P)|\geq 4n-(n-1)=3n+1$. By Lemma \ref{L2}, $(G[U\less V(P)],\tau)$ must contain    a red edge, say $uv$.       But  then $(G,\tau)$ has a red $C_{2n+1}$  with edge set  $E(P)\cup\{xu_1, xu, uv, vy, yu_j\}$, contrary to Claim~\ref{noC2n+1}. Thus $H$ has no  path $P$ on $2n-3$ vertices     with both ends in  $V(C)$.  It follows that $|W_1|\le n-3$ and $|W_2|\le n-3$.  \medskip

We next show  that  $|W_0|\ge1$. Suppose $W_0=\emptyset$. Then $|W_1|+
|W_2|=|W|=n$. We may assume that  $|W_1|\le |W_2|$.  Then $  |W_1|\ge3$ because $|W_2|\le n-3$.  Recall that  $|U_1|=n$
and $|V_p|\le n-1$,   we see that $(G[U_1],\tau)$ contains an 
edge $e$ that is colored   red or    blue. We may assume that $e= u_1u_{2j+1}$  for some $j\in [n-1]$.  Then  $e$ must be red, else  let  $v\in W_1$,   we obtain a blue   $C_{2n+1}$ with vertices $v, u_{2n},u_{2n-1},\ldots,u_{2j+1},u_1,u_2,\ldots, u_{2j}$ in order, contrary to Claim~\ref{noC2n+1}.  
Let $v_1\in U_1\less\{u_1, u_{2j+1}\}$ and $v_2 \in U_2\less u_2$; let 
$Q_1$ be a   $(u_1, v_1)$-path on $2|W_1|+1$ vertices using  edges between $W_1$ and
$U_1\less u_{2j+1}$ and  $Q_2$ be a  $(u_2, v_2)$-path on $2|W_2|-3$
vertices  using  edges between $W_2$ and
 $U_2$. But then $(G,\tau)$ contains a red $C_{2n+1}$  with edge set  $E(Q_1)\cup E(Q_2)\cup\{xu_{2j+1}, u_{2j+1}u_1,  
v_1y, yv_2, u_2x\}$, contrary to Claim~\ref{noC2n+1}.  This proves that    $|W_0|\ge 1$, as claimed. \medskip

Finally, let $H'$ be obtained from $H$ by adding a new vertex $w$ adjacent to all vertices in $W$. Then $H'$ is a bipartite graph with bipartition $(V(C)\cup\{w\}, W)$. Note that for each vertex $v\in W$, $d_{H'}(v)\ge (|V(C)\cup\{w\}|+1)/2$;    for each vertex $v\in W_0\ne \es$, $d_{H'}(v)> (|V(C)\cup\{w\}|+1)/2$. By Lemma \ref{Hal2} applied to $H'$ with $M:=V(C)\cup\{w\}$, $N:=W$ and $\ell:=n$, we see that $H'$ contains a cycle $C_{2n}$  and thus $H$ has a path $P$ on $2n-3$ vertices with both ends on $C$, contrary to the fact that $H$ has no  path  on $2n-3$ vertices     with both ends in  $V(C)$. 
\medskip

It remains to consider the case   $|R|\leq 4n-1$.  In this case,  
\[
 |B|=|V(G)\less (X\cup  V_p\cup  R)|\ge n\cdot 2^k+1-(k-1)(n-1)-(n-1)-(4n-1) \ge n+5.
 \]
By Lemma \ref{L},   $|X_1\cup V_p|\leq n$ and $|X_2\cup V_p|\leq n$. Since $X_1=\emptyset$ or $X_2=\emptyset$, we see that $|X_1\cup X_2\cup V_p|\le n$.  By  Claim \ref{Odd_empty},   $|X\cup V_p|=|X_3|+\cdots+|X_k|+|X_1\cup X_2\cup V_p|\leq (k-2)(n-1)+n$. Thus 
\[
|B|+|R|= |G|-|X\cup V_p|\geq  n\cdot 2^k+1-(k-2)(n-1)-n\ge 6n+2,
\]
which implies $|B|\ge 2n+3$. Then $|R|\ge|B|\ge 2n+3$.  By Lemma~\ref{L2}, $(G[R],\tau)$ must have a red edge, say $x_ry_r$;  $(G[B],\tau)$  must have a blue edge, say $x_by_b$.  Let  $W\subseteq R\less \{x_r, y_r\} $ be such that $|W|=2n+1$.  Then all the edges between $B$ and $W$ are colored     red or blue.   Let 
\[
\begin{split}
B_r: &=\{v\in B\less \{x_b, y_b\}  \mid v \text{ is red-adjacent to at least  } n+1 \text{  vertices in  }   W\} \text{ and } \\
B_b: &=\{v\in B\less \{x_b, y_b\}   \mid v \text{ is blue-adjacent to at least  } n+1 \text{  vertices in  }   W\}. \\
\end{split}
\]
Then  $|B_r|\geq n+1$ or $|B_b|\geq n+1$  because $|B\less \{x_b, y_b\}|\geq 2n+1$. We may assume that $|B_j|\geq n+1$ for some color    $j\in\{r, b\}$.   Let $H_j$ be the bipartite subgraph of $G$ with bipartition  $(W, B_j)$ such that $E(H_j)$ consists of  all   edges in color $j$  between $W $ and $B_j$ in $(G,\tau)$. For each vertex $v\in B_j$, $d_{H_j}(v)\ge (|W|+1)/2$.    By Lemma \ref{Hal2} applied to $H_j $ with $M:=W$, $N:=B_j$ and $\ell:=n$, we see that $H_j $ has a    path $P$ in color $j$ on $2n-3$ vertices with  both ends in $ W$ when $j=r$ and both ends in $B_j$ when $j=b$.   Let  $v_1, v_2$ be the two ends of $P$. But then $(G,\tau)$ contains  a monochromatic copy of  $C_{2n+1}$ in color $j$ with edge set  $E(P)\cup  \{v_1x, xx_j,x_jy_j, y_jy, yv_2\} $,  contrary to Claim~\ref{noC2n+1}.  \qed
 \medskip

 \begin{claim}\label{Odd_Vp2}
 $|V_{p-2}| \le n$ if $p\ge 3$.
\end{claim}
\pf  Suppose $|V_{p-2}|\ge n+1$. Then $n+1 \le |V_{p-2}| \le |V_{p-1}| \le |V_p|$.  Let $B_1$, $B_2$, $B_3$ be a permutation of $V_{p-2}$, $V_{p-1}$, $V_p$.  Since $(G, \tau)$ has no rainbow triangle, we may assume that     $B_2$ is, say blue-complete,  to $B_1 \cup B_3$ in $(G, \tau)$. Then $B_1$ must be  red-complete to $B_3$ in $(G, \tau)$, else we obtain a blue    $C_{2n+1}$, contrary to Claim~\ref{noC2n+1}.     By Lemma~\ref{L}, no vertex in $X$ is either red- or blue-complete to $V(G)\less X$. Thus $X_1=X_2=\emptyset$ and $|X| \le (k-2)(n-1)$. Let $A:=V(G)\less (B_1\cup B_2\cup B_3\cup X)$.   By Lemma~\ref{L},  no vertex in $A$ is red-complete to $B_1\cup B_3$ in $(G, \tau)$,  and no vertex in $A$ is blue-complete to $B_1\cup B_2$ or $B_2\cup B_3$ in $(G, \tau)$. This implies that 
$A$ must be  red-complete to $B_2$ in $(G, \tau)$.    Suppose that $(G[A],\tau)$ has a blue edge, say, $uv$.  We may assume that $u$ is blue-complete to $B_1$. But then we   obtain a blue   $C_{2n+1}$, because $v$ is blue-complete to $B_1$ or $B_3$. 
This proves that  $(G[A],\tau)$ has no blue edges. By Lemma~\ref{L}, $(G[B_2],\tau)$ has no blue edges, and neither $(G[B_1],\tau)$ nor $(G[B_3],\tau)$ has red or blue edges. By minimality of $k$, 
$ |B_1|\le n \cdot 2^{k-2}$ and $|B_3|\le   n \cdot 2^{k-2}$. Let $\tau'$ be obtained from $\tau$ by recoloring all the blue edges in $G[X]$ red.  Then $\tau'$  is  a Gallai $(k-1)$-coloring of $G[B_2\cup A\cup X]$. Moreover, $(G[B_2\cup A\cup X],\tau')$ has   no blue edges and no monochromatic copy of $\Theta_{2n+1}$ when $n\ge3$ and no monochromatic $C_5$ when $n=2$. By minimality of $k$, $|B_2\cup A\cup X|\le n\cdot 2^{k-1}$.    But then 
\[
 |G|=|B_1 |+|B_2\cup A\cup X|+|B_3 | \le n \cdot 2^{k-2}+n\cdot 2^{k-1} +n\cdot 2^{k-2} <n\cdot 2^{k}+1,
 \] 
 for all $k\ge 3$ and $n\ge 2$, a contradiction.  
   \qed

\begin{claim}\label{Odd_BR}
 $|R| \ge n+1$ and $|B| \ge n+ 1$.  
\end{claim}
\pf  Suppose  $|R| \le n $ or $|B| \le n$, say the latter.  By Claim \ref{Odd_Vp} and \ref{Odd_B_R}, $|V_p|\ge n+1$ and $|R| \ge n+2$.  By Lemma \ref{L} applied to $V_p$ and $R$,   we see that $X_1 = \emptyset$ and neither $(G[V_p], \tau )$
nor $(G[R], \tau)$ has red edges.   Let  $i\in[k]$ be a color, say green,   that is neither  red nor blue. Let $\tau'$ be obtained from $\tau$ by recoloring all the red edges in $(G[B\cup (X\backslash X_i)], \tau)$ by color green; and recoloring  all the red edges in $(G[X_i], \tau)$ by color blue.  Then neither $(G[R\cup X_i],\tau')$ nor $(G[B\cup V_p\cup (X\less X_i)],\tau')$ has   red edges.   It follows that  $\tau'$ is  a Gallai $(k-1)$-coloring of  both $G[R\cup X_i]$  and  $G[B\cup V_p\cup (X\less X_i)]$, and neither $(G[R\cup X_i],\tau')$ nor $(G[B\cup V_p\cup (X\less X_i)],\tau')$ has   a monochromatic copy of $\Theta_{2n+1}$ when $n\ge3$ and   monochromatic $C_5$ when $n=2$.  By minimality of $k$, $|R\cup X_i|\le n\cdot 2^{k-1}$ and $|B\cup V_p\cup (X\less X_i)|\le n\cdot 2^{k-1}$.  Then 
\[
|G|=|R\cup X_i|+|B\cup V_p\cup (X\less X_i)|\leq  n\cdot 2^{k-1 }+n\cdot 2^{k-1 }<n\cdot 2^k+1, 
\]
for all $k\ge 3$ and $n\ge 2$, a contradiction.  
 \qed\medskip

By Claim~\ref{Odd_BR},  $|R| \ge n+1$ and $|B| \ge n+ 1$.  By Claim \ref{Odd_Vp}, $|V_p|\ge n+1$. By Lemma \ref{L} applied to $V_p$ and $R$, and $V_p$ and $B$, respectively, we have $|X_1|=|X_2|=0$;     $(G[R],\tau)$ has no red edges; $(G[B],\tau)$ has no blue  edges;  and $(G[V_p],\tau)$ has neither red nor blue edges.  By Claim~\ref{Odd_X_i bound},    $|X| \le  (k-2)(n-1) $.   By minimality of $k$,   $|V_p| \le n\cdot 2^{k-2}$.   We may further assume that $V_{p-1}\subseteq B$.  By Claim~\ref{Odd_Vp2},   $|V_{p-2}|\le n$ if $p\ge 3$.   By Lemma \ref{L2},    $|B \less V_{p-1}| \le 2n$ and $|R| \le 2n$.    If $k\ge4$ or $|X|=0$, then
\[
\begin{split}
|G| &= |V_p| + |V_{p-1}| + |B \less V_{p-1}| + |R| + |X| \\
&\le \begin{cases}  
n\cdot 2^{k-2} +n\cdot 2^{k-2}+ 2n  + 2n +(k-2)(n-1)  &  \text{if } k\ge4 \\
n\cdot 2^{k-2} +n\cdot 2^{k-2}+ 2n  + 2n      &  \text{if } k=3, |X|=0\\
\end{cases}\\		
	&< n\cdot 2^{k}+1, 
\end{split}
\]
  a contradiction.  Thus $k=3$ and $X\ne \emptyset$. We may assume that $X$ is green-complete to $V(G)\less X$. Then all the edges of $G[V_p]$ are colored green. Thus  $|X\cup V_p|\le 2n$, else we obtain a green copy of $\Theta_{2n+1}$.  It follows that  $|V_{p-1}|\le |V_p|\le  2n-1$. But then  
  \[
 |G| = |X\cup V_p| + |V_{p-1}| + |B \less V_{p-1}| + |R| \le 2n+(2n-1)+2n+2n<8n+1,
  \]
  for all $n\ge2$,  a contradiction.  \medskip
  
     This completes the proof of Theorem~\ref{oddchorded}. \qed

\section*{Acknowledgements}

The authors would like to thank Christian Bosse and Jingmei Zhang for their helpful discussion.  

\end{document}